\documentclass[12pt]{article}
\usepackage{mathpazo}
\usepackage{amsfonts,amssymb,amsmath}
\usepackage{enumerate}
\usepackage{epigraph}
\usepackage{color,soul}     
\usepackage[normalem]{ulem}  

\newtheorem{proposition}{Proposition}

\newtheorem{lemma}{Lemma}

\newtheorem{consequence}{Consequence}
\newtheorem{theorem}{Theorem}
\newtheorem{definition}{Definition}

\newtheorem{question}{Question}
\newtheorem{thesis}{Thesis}
\newenvironment{proof}[1][Proof]{\noindent\textbf{#1.} }{\ \rule{0.5em}{0.5em}}

\begin{document}

\title{\textbf{Sets and Probability}\thanks{Versions of this paper have been presented a Bristol--Leuven workshop on Logic and Philosophy of Science (2015), at the Philosophy Department of the Universidade Federal do Rio
Grande do Norte (2015), the Fourth Reasoning Conference in Manchester (2015), the Philosophy of Mathematics Seminar in Oxford (2014), and at the Philosophy Departmental Research Seminar in Aberdeen (2014). We are grateful to the audiences for helpful comments, questions, and suggestions. In this respect we are especially indebted to Philip Welch, George Wilmers, and Sylvia Wenmackers.}}
\author{Hazel Brickhill and Leon Horsten}
\date{\today}
\maketitle

\begin{abstract}
\noindent In this article the idea of random variables over the set theoretic universe is investigated. We explore what it can mean for a random set to have a specific probability of belonging to an antecedently given class of sets.

\end{abstract}


\section{Introduction}\label{intro}


Probabilistic notions have been applied to mathematical objects and notions. For instance, probabilistic concepts have been applied in the theory of random graphs \cite{AS}. The aim of this article is to apply a notion of probability to the mathematical universe as a whole. More in particular, 
we wish to explicate \emph{what it could mean for a property $A$ of sets to have a probability of being true of a set $y$ in the set theoretic universe $V$.} 
Properties are identified with their extensions, so that $A$ ranges over all proper and improper classes in $V$. 

The aim is to develop a theory of the probability of events of the form $A(\tau)$, where $A$ is a class and the variable $\tau$ is a \emph{random variable}. The \emph{state space} of the random variables is of course $V$. The \emph{outcome space} of the random variables has to be at least as large as $V$ because there must be enough states for a random variable to take each set as a possible value. On the other hand, there is no need for it to be larger than $V$. Therefore the outcome space is simply identified with $V$.


Without invoking fixed set of postulates, intuitions about probability have occasionally been used in set theory,  for instance to motivate new basic principles \cite{F}. However, such attempts are mostly regarded as unsuccessful \cite{H}. In the light of this it is natural to wonder what we should require from probability functions associated with random variables on $V$.

Surely it would be unreasonable to insist on there being one \emph{unique} correct probability function that yields the probability of a random variable taking a value in a given class of sets. On the other hand, for our functions to have any hope of meriting the label \emph{probability function}, they have to satisfy Kolmogorov's conditions for being a \emph{finitely additive probability function}.

From the outset we impose additional constraints on the class of probability functions that we are interested in:\footnote{For a discussion of these constraints in the context of non-Archimedean probability theory, see \cite{HWB}.}
\begin{enumerate}
\item \textbf{Totality.} The probability functions are defined on all classes.
\item \textbf{Uniformity.} All singleton events are given the same probability.
\item \textbf{Regularity.} All singleton events are given non-zero probability.
\end{enumerate}
All this means, for familiar reasons, that the sought-for probability functions cannot be Kolmogorov probability functions. Given our insistence on finite additivity, this means that the probability functions will be non-Archimedean. They will not satisfy $\sigma$-additivity, but they will instead satisfy a generalised infinite additivity rule.

In mathematics today, the term `probability' has become virtually \emph{synonymous} with `function that satisfies the Kolmogorov axioms (\emph{including} $\sigma$-additivity)'. If you see  matters this way, then you will will be loath to dignify the functions constructed in this paper by the term `probability function'. Nonetheless, you  may ask the question whether a fine-grained quantitative theory of possibility, with which the degree of possibility of properties can quantitatively be compared, can be constructed. This is what is investigated in the present article. So, if you prefer, you can call the theory constructed in this paper a \emph{quantitative theory of possibility.} You are then advised to replace all occurrences of `(non-standard) probability function' by `quantitative possibility function'.

The project in which we are engaging in this article is related to the work in \cite{BDF}. The aim of the latter article is to construct a theory of sizes for mathematical universes inspired by the \emph{Euclidean principle} that the size of the whole is larger than the sizes of its proper parts. 
Now there is of course a familiar theory of size---Cantor's theory of cardinality,---which does not satisfy this Euclidean principle. So Benci and his co-authors propose their Euclidean theory of size as a rival to Cantor's theory.

We, on the other hand, fully accept Cantor's theory of cardinality. None\-theless, the probability functions that will be constructed satisfy the Euclidean principle that the probability of an event is \emph{strictly greater} than the probability of each of its sub-events. Moreover, the mathematical techniques for generating them are closely related to the techniques that are used in \cite{BDF}.

What we shall mean by `mathematical universe' is not the same as what is meant \cite{BDF} by the term. The authors of \cite{BDF} impose mainly algebraic constraints on what counts as a mathematical universe \cite[Introduction]{BDF}. We, in contrast, take the term `mathematical universe' in the set theoretical sense. Naively, you may take there to be one preferred set theoretic universe: $V$. But if you are uncomfortable with taking $V$ as given, then you might want to take a mathematical universe to be a rank $V_{\alpha}$ that constitutes a model of most or perhaps even all of the standard principles of set theory. Indeed, we will see that for random variables defined on \emph{any} large set $S$, the general idea of equipping them with a probability function will be the same as that for random variables on $V$. 

We will discuss two ways of generating non-Archimedean probability functions for random variables on $V$. In section \ref{snapshot} a simple way of generating such probability functions (the \emph{finite snapshot approach}) will be described. In section \ref{constraints} we go on to discuss how global properties of these probability functions can be made to hold by imposing constraints on the process of generating such functions. In section \ref{bootstrapping}, a theoretically more satisfying but also more complicated way of generating non-Archimedean probability functions for random variables on $V$ is discussed (the \emph{bootstrapping} method).


\section{The finite snapshot approach}\label{snapshot}

A random variable $\tau$ on $V$ is a function from states to the outcome space, i.e., an element of ${^V}V$. So there are \emph{many} random variables on $V$. The aim is to associate a notion of probability with elements of ${^V}V$ that meet the minimal constraints (totality, uniformity and regularity) that were described in section \ref{intro}.

In fact, we want to give precise meaning to conditional probability\index{conditional probability} statements of the form
$$ {\mathsf{Pr}(\sigma \in A \mid \tau \in B)       },$$
where $\sigma, \tau \in {^V}V$ and $A,B \subseteq V$. But we will see that it will be sufficient for our purposes to give meaning to \emph{unconditional} probability statements of the form $ {\mathsf{Pr}(\sigma \in A )}.$ So our fundamental problem amounts to giving meaning to expressions of the form $ \mathsf{Pr}(\sigma \in A ).$ Such probability measures will be determined by a choice of a fine ultrafilter\index{ultrafilter} on the collection $[V]^{<\omega}$ of \emph{finite subsets} of the state space\index{state space}.\footnote{What follows is an adaptation of the approach of \cite[section 2]{BrickhillHorsten2018}.}

The starting point is a fine ultrafilter\ $\mathcal{U}$ on $[V]^{< \omega}$. 
This fine ultrafilter\index{ultrafilter!fine} $\mathcal{U}$ defines a non-Archimedean field $\mathcal{F}_{\mathcal{U}}$ in the following way.

For any two functions $f,g: [V]^{< \omega } \rightarrow \mathbb{Q}$ we define:

\begin{definition}
$$f \approx_{\mathcal{U}} g \equiv \{ T \in  [V]^{< \omega } : f(T) = g(T)  \} \in \mathcal{U}.$$
\end{definition}
\noindent In words: two functions are identified if they coincide on ultrafilter-many states.

The relation $ \approx_{\mathcal{U}}$ is an equivalence relation, so we can take equivalence classes  for which we then have
$$ [f]_{\mathcal{U}} = [g]_{\mathcal{U}} \Leftrightarrow   f \approx_{\mathcal{U}} g   .$$
Moreover, it is again a routine exercise to verify that the $[f]_{\mathcal{U}}$'s form a hyper-rational field $\mathcal{F}_{\mathcal{U}}$.

Now suppose $A \subseteq V$ and $\theta \in {^V}V$. Then we define the function $f_{\theta \in A}: [V]^{< \omega } \rightarrow \mathbb{Q}$ as follows:
\begin{definition}
For every $T \in  [V]^{< \omega } : $
$$ f_{\theta \in A}(T) \equiv \frac{\vert  \{ s\in T: \theta (s) \in A     \} \rvert}{\lvert  T  \rvert}     .$$
\end{definition}
\noindent In words: for every \emph{finite} set of states $T$, $f_{\theta \in A}(T)$ is the ratio between the number of states $s$ in $T$ for which $\theta (s) \in A$ and the number of states in $T$. In this sense, $f_{\theta \in A}(T)$ is the probability of $\theta \in A$ \emph{on a finite snapshot of states}.

Similarly, we define the function $f_{\theta \in A \wedge \nu\in B}$ as follows:
\begin{definition}
For every $T \in  [V]^{< \omega } : $
$$ f_{\theta \in A \wedge \nu \in B}(T) \equiv \frac{\vert  \{ s\in T: \theta (s) \in A  \textrm{ and } \nu (s) \in B   \} \rvert}{\lvert  T  \rvert}     .$$
\end{definition}

Now we are ready to define the probability of  $\theta \in A$, relative to a fine (and therefore free) ultrafilter $\mathcal{U}$ on $[V]^{< \omega}$:

\begin{definition}
$$ \mathsf{Pr}_{\mathcal{U}} (\theta \in A) \equiv     [ f_{\theta \in A} ]_{\mathcal{U}} .  $$
\end{definition}
\noindent Similarly, we define $ \mathsf{Pr}_{\mathcal{U}} (\theta \in A \wedge \nu \in B)$ as $[ f_{\theta \in A \wedge \nu \in B} ]_{\mathcal{U}}$.
Thus we have constructed a probability function $ \mathsf{Pr}_{\mathcal{U}}$ that takes its values in the hyper-rational field $\mathcal{F}_{\mathcal{U}}$. Such probability functions are sometimes called NAP functions.

Conditional probability can then be expressed in terms of unconditional probability:
\begin{definition}
$$  \mathsf{Pr}_{\mathcal{U}} (\theta \in A \mid \nu \in B)   \equiv  \frac{\mathsf{Pr}_{\mathcal{U}} (\theta \in A \wedge \nu \in B)}{\mathsf{Pr}_{\mathcal{U}} (\nu \in B)}  .$$
\end{definition}


\section{Constraints}\label{constraints}

From section \ref{intro} we know that the aim is not to arrive at a \emph{unique} (correct) probability function on $V$. But we did insist from the outset on our probability functions satisfying three global constraints: totality, uniformity, and regularity. It will be shown that these properties are always guaranteed to hold.

There are further global conditions on probability functions on $V$ that seem reasonable to require, and that are not guaranteed to hold without further work. These global constraints will be explored. We will show that many of them can be forced to hold by imposing constraints on the ultrafilters from which the probability functions are generated.

\subsection{Elementary properties}

The definition of $ \mathsf{Pr}_{\mathcal{U}}$ is relative to an initial choice of the fine ultrafilter $\mathcal{U}$. The properties of $ \mathsf{Pr}_{\mathcal{U}}$ depend on $\mathcal{U}$. Nonetheless, certain basic properties of $ \mathsf{Pr}_{\mathcal{U}}$ can be easily seen to hold regardless of which fine ultrafilter $\mathcal{U}$ is chosen:

\begin{proposition}\label{basic properties} \textrm{}
\begin{enumerate}
\item $ \mathsf{Pr}_{\mathcal{U}}$ is a finitely additive probability function;
\item  $ \mathsf{Pr}_{\mathcal{U}}$ is Euclidean.
\end{enumerate}

\begin{proof}
Easy.
\end{proof}
\end{proposition}

Now we define the notion of a \emph{diagonal random variable}:
\begin{definition}
A random variable $\theta$ is said to be a \emph{diagonal} random variable if for any set $x$, there is exactly one element $u$ of the state space such that $\theta (u) = x$.
\end{definition}
\noindent In words: a diagonal random variable is a random variable that takes every value exactly once.

Using this notion, we define the notions of \emph{regularity} and \emph{uniformity}:
\begin{definition}[regularity]
A probability function $ \mathsf{Pr}_{\mathcal{U}}$ is regular if for every diagonal random variable $\theta$ and for every $x\in V, \mathsf{Pr}_{\mathcal{U}}(\theta = x) >0$.
\end{definition}

\begin{definition}[uniformity]
A probability function $ \mathsf{Pr}_{\mathcal{U}}$ is uniform if for every diagonal random variable $\theta$ and for  all $x,y \in V:$ $$\mathsf{Pr}_{\mathcal{U}}(\theta = x) =  \mathsf{Pr}_{\mathcal{U}}(\theta = y) . $$
\end{definition}

\begin{proposition} \textrm{}
For every fine ultrafilter $\mathcal{U}$: 
\begin{enumerate}
\item $ \mathsf{Pr}_{\mathcal{U}}$ is regular;
\item $ \mathsf{Pr}_{\mathcal{U}}$ is uniform.
\end{enumerate}

\begin{proof}
These properties are proved as propositions 2.5 and 2.6 in \cite[p.~525--526]{BrickhillHorsten2018}.
\end{proof}
\end{proposition}

The Euclidean property is formally defined as follows:
\begin{definition}[Euclidean]
A probability function $ \mathsf{Pr}_{\mathcal{U}}$ is Euclidean if for every diagonal random variable $\theta$ and all $A,B \subseteq V$:
$$ A \subsetneq B \Rightarrow  \mathsf{Pr}_{\mathcal{U}} (\theta \in A) < \mathsf{Pr}_{\mathcal{U}} (\theta  \in B).$$
\end{definition}

Then we have:
\begin{proposition}
For every fine ultrafilter $\mathcal{U}$, the probability function $ \mathsf{Pr}_{\mathcal{U}}$ is Euclidean.

\begin{proof}
By finite additivity and regularity.
\end{proof}
\end{proposition}

Now we turn to infinite additivity.
Countable additivity\index{countable additivity} means that the probability of the union of a countable family of disjoint sets is the \emph{infinite sum}\index{infinite sum} of the probabilities of the elements of the family, where the notion of infinite sum is spelled out in terms of the classical notion of limit\index{limit!classical}.
In the present setting, the probability $Pr_{\mathcal{U}}$ of the union of \emph{any} family of disjoint sets is also the infinite sum of the probabilities of the elements of the family \cite[section 3.4]{BHW}. But now the notion of infinite sum is spelled out in terms of the generalised notion of limit\index{limit!generalised} based on the ultrafilter\index{ultrafilter} $\mathcal{U}$. More precisely, the new notion of infinite sum is defined as follows. Suppose we are given a family $\{ q_i: i \in \mathbb{N} \}$ of rational numbers, and $I\subseteq \mathbb{N}$. Then consider the function $f: [\mathbb{N}]^{< \omega} \rightarrow \mathbb{Q}$ given by
$$ f(T) = \sum_{i \in I \cap T} q_i .$$
This function can be seen as giving the value of the infinite sum on all \emph{finite parts} (``snapshots'') of the index set. So we identify the infinite sum of the family $\{ q_i: i \in I \}$ of rational numbers with the generalised limit\index{limit!generalised} of $f$ according to the ultrafilter\index{ultrafilter} $\mathcal{U}$:
\begin{definition}
$$  {\sum_{i\in I}}^* q_i \equiv [f]_{\mathcal{U}}.   $$
\end{definition}

Using this notion of infinite sum\index{infinite sum}, we can express the probability of the union of a disjoint family of sets as the sum of the probabilities of the members of that family:
\begin{proposition}\label{perfectly additive}\index{probability function!perfectly additive}
If $A = \bigcup_{i\in I}A_i$, with $A_i \cap A_j = \emptyset$ for all $i,j \in I$, then for every random variable $\tau$:
$$\mathsf{Pr}_{\mathcal{U}}(\tau \in A) =  {\sum_{i\in I}}^*\mathsf{Pr}_{\mathcal{U}}(\tau \in A_i).$$
\end{proposition}
\noindent In sum, $\mathsf{Pr}_{\mathcal{U}}$ has a natural infinite additivity property that is sometimes called \emph{perfect additivity}.

\begin{proposition}
For every fine ultrafilter $\mathcal{U}$, the probability function $\mathsf{Pr}_{\mathcal{U}}$ is perfectly additive.

\begin{proof}
This proposition is proved as proposition 8 in \cite[p.~132--133]{BHW}.
\end{proof}
\end{proposition}


\subsection{Symmetry principles}

From now on, the symbol $\theta$ will be used to refer to some arbitrary \emph{diagonal} random variable. When it is not assumed that the random variable in question is diagonal, we will write $\tau$.

The Euclidean-ness of  $ \mathsf{Pr}_{\mathcal{U}}$  has implications for \emph{symmetry principles}\index{symmetry principle}. As a rule of thumb, one can say that \emph{symmetry principles fail}.\index{symmetry principle}\footnote{See  \cite{BDF}, \cite{BHW}, \cite{HWB}.} 

\begin{proposition}\label{invar}
For every fine ultraflter $\mathcal{U}$, the probability function $Pr_{\mathcal{U}}$ is not invariant under all permutations of $V$.

\begin{proof}
We concentrate on $\mathbb{N}$ as it is canonically represented in $V$ (by means of the Zermelo ordinals, for instance). Define a permutation $\pi$ of $V$ as follows:
\begin{itemize}
\item $\pi(x) = x$ for $x \in V \setminus \mathbb{N}$; Otherwise:
\item $\pi(x) = x+2$ for $x$ even;
\item $\pi (1) = 0$;
\item $\pi ( x) = x -2$ for $x$ odd and $>1$.
\end{itemize}
Let $A \equiv \{0,2,4, \ldots  \}$, and let $\theta$ be a diagonal random variable. Then $\pi(A) \subsetneq A$. Therefore, by the Euclidean principle, $ \mathsf{Pr}_{\mathcal{U}}(\theta \in \pi(A) ) < \mathsf{Pr}_{\mathcal{U}}( \theta \in A)  .$
\end{proof}
\end{proposition}
\noindent This of course entails that there are diagonal random variables $\theta, \theta'$ such that for some $A\subseteq V$, $$ \mathsf{Pr}_{\mathcal{U}}( \theta \in A ) \neq \mathsf{Pr}_{\mathcal{U}}(\theta' \in A ).$$

One popular global constraint on probability measures is \emph{translation-invariance}\index{translation-invariance}. The Lebesgue\index{Lebesgue, Henri} measure\index{Lebesgue measure}\index{Lebesgue, Henri} has this property, and Banach limits\index{Banach limit}\index{Banach, Stefan} seem to occupy a privileged position in the class of generalised limits\index{limit!generalised} at least in part because they are translation-invariant\index{translation-invariance}. In our context, translation-invariance\index{translation-invariance} does not make obvious sense. For a random class $A$, it is not clear what `$A + \alpha$' (where $\alpha$ is a number) \emph{means}. But a clear interpretation of `adding an ordinal number' can of course be given if $A$ is a collection of ordinals:
\begin{definition}
For $A$ any collection of ordinals:
$$ A\oplus \alpha \equiv \{ \beta : \exists \gamma \in A \textrm{ such that } \beta = \gamma + \alpha  \}  .$$
\end{definition}
\noindent Then for $A$ to be translation-invariant\index{translation-invariance} means that for all ordinals $\alpha$ and for every $\theta$, $$ \mathsf{Pr}_{\mathcal{U}}(\theta \in A) = \mathsf{Pr}_{\mathcal{U}}(\theta \in A\oplus \alpha).$$

However, even if we consider non-Archimedean measures\index{probability!non-Archimedean} (of the kind that we have been describing) on ordinals, translation-invariance conflicts with the Euclidean Property of our generalised probability functions. In particular, there is no $\mathrm{NAP}$ probability function $\mathsf{Pr}_{\mathcal{U}}$ on any infinite cardinal $\kappa$ such that there is even one ordinal $\alpha$ with $0 < \alpha < \kappa$ and $$\mathsf{Pr}_{\mathcal{U}}(\theta \in \kappa) = \mathsf{Pr}_{\mathcal{U}}(\theta \in \kappa \oplus \alpha).$$ The reason is simple. We have $\kappa \oplus \alpha = \kappa \backslash \alpha \subsetneq \kappa,$ so if we had $\mathsf{Pr}_{\mathcal{U}}(\theta \in \kappa) = \mathsf{Pr}_{\mathcal{U}}(\theta \in \kappa \oplus \alpha),$ then we would contradict the Euclidean principle.

As this example shows, such translations arenÕt necessarily one to one so we may not want full invariance in general.
In \cite[section 1.3]{BDF}, Benci, Forti,\index{Forti, Marco} and Di Nasso\index{Di Nasso, Mauro} explore a \emph{restricted} notion of translation-invariance\index{translation-invariance} of $\mathrm{NAP}$-like measures on ordinals. We do not pursue this theme further here, but only pause to note that there are other reasonable-looking principles that are hard to satisfy. In the context of their theory of \emph{numerosities}, Benci,\index{Benci, Vieri} Forti\index{Forti, Marco}, and Di Nasso\index{Di Nasso, Mauro} consider a principle that in the present context would take the following form:
\begin{definition}[Difference Principle]\index{difference principle} \textrm{}

\begin{multline*}  \forall A, B \in V: \mathsf{Pr}_{\mathcal{U}}(\theta \in A) < \mathsf{Pr}_{\mathcal{U}}(\theta \in B) \Rightarrow \\ \exists C \in V:  \mathsf{Pr}_{\mathcal{U}}(\theta \in B) = \mathsf{Pr}_{\mathcal{U}}(\theta \in A) + \mathsf{Pr}_{\mathcal{U}}(\theta \in C) .  \end{multline*}
\end{definition} On countable sample spaces, the difference principle\index{difference principle} can be made to hold by building $\mathsf{Pr}_{\mathcal{U}}$ from a \emph{selective} ultrafilter\index{ultrafilter!selective} \cite{BD}. But the existence of selective ultrafilters\index{ultrafilter!selective} is independent of $\textrm{ZFC}$. As far as we know, it is an open whether the difference principle\index{difference principle} can be consistently made to hold for $\textrm{NAP}$ probability functions on uncountable sample spaces.


\subsection{Probability and cardinality}

In this (sub-)section we investigate the relation between our notion of generalised probability on the one hand, and the familiar notion of cardinality on the other hand.

\subsubsection{Hume's principle for probability}

One might naively wonder whether the following probabilistic analogue of Hume's Principle for cardinality can hold:

\begin{definition}[Hume's principle for probability]\label{Hume}
For all $A,B \in V$: $$\left \vert A \right \vert = \left \vert B \right \vert \Rightarrow \mathsf{Pr}_{\mathcal{U} }(\tau \in A) = \mathsf{Pr}_{\mathcal{U} }(\tau \in B).$$
\end{definition}

But the probability functions $\mathsf{Pr}_{\mathcal{U}}$ that we have been considering cannot satisfy Hume's principle for probability, as its failure is an immediate consequence of Proposition \ref{invar}: invariance under permutations and Hume's principle for probability are mathematically equivalent. 
However, this was only to be expected. After all, we do not expect Kolmogorov probability (on infinite spaces) to satisfy any such principle.


\subsubsection{Superregularity}

The hyper-rational field $\mathcal{F}_{\mathcal{U}}$ in which the probability functions $\mathsf{Pr}_{\mathcal{U}}$  take their values contain \emph{infinitesimal numbers}\index{infinitesimal}---this is what makes it non-Archi\-me\-dean. We will write $\mathsf{Pr}_{\mathcal{U}}(\sigma \in A)\approx 0$ if $\mathsf{Pr}_{\mathcal{U}}(\sigma \in A) < n^{-1}$ for each $n \in \mathbb{N}$. And we will write $\mathsf{Pr}_{\mathcal{U}}(\sigma \in A) \ll \mathsf{Pr}_{\mathcal{U}}(\tau \in B)$  if $$\frac{\mathsf{Pr}_{\mathcal{U}}(\sigma \in A)}{\mathsf{Pr}_{\mathcal{U}}(\tau \in B)} \approx 0.$$

We have seen that $\mathsf{Pr}_{\mathcal{U}}$ cannot satisfy Hume's principle for probability. But, at least at first sight, it seems that it would be reasonable to demand:
 $$ \left \vert A \right \vert < \left \vert B \right \vert \Rightarrow \mathsf{Pr}_{\mathcal{U}}(\delta \in A) < \mathsf{Pr}_{\mathcal{U}}(\delta \in B) .$$ Indeed, if in addition $\left \vert B \right \vert \geq \omega$, then we might even expect $$\left \vert A \right \vert < \left \vert B \right \vert \Rightarrow \mathsf{Pr}_{\mathcal{U}}(\sigma \in A) \ll \mathsf{Pr}_{\mathcal{U}}(\sigma \in B).$$ Further, this may be expected to hold if $B$ is a proper class but $A$ is a set .\index{proper class} The result is a size constraint which is a strengthening of the requirement of regularity:

\begin{definition}[Superregularity]\label{superregularity}
$$ \omega \leq \left \vert A \right \vert < \left \vert B \right \vert \leq \left \vert V \right \vert   \Rightarrow \mathsf{Pr}_{\mathcal{U}}(\theta \in A) \ll \mathsf{Pr}_{\mathcal{U}}(\theta \in B) .$$ 
\end{definition}

Note that if $A$ is finite and $B$ is infinite then the consequent holds automatically.

By a suitable restriction on admissible ultrafilters\index{ultrafilter} $\mathcal{U}$, superregularity\index{superregularity} can indeed be made to hold:
\begin{theorem}\label{superregtheorem}
There are fine ultrafilters\index{ultrafilter} $\mathcal{U}$ such that $ \mathsf{Pr}_{\mathcal{U}}$ is superregular.

\begin{proof}

If $A,B\in V$ such that $\omega \leq \left \vert A \right \vert < \left \vert B \right \vert $ are given, then we have $\mathsf{Pr}_{\mathcal{U}}(\theta \in A) \ll \mathsf{Pr}_{\mathcal{U}}(\theta \in B)$ if and only if for each $n \in \mathbb{N}$,  $$ \{ D \in [V]^{< \omega} : \frac{\mathsf{Pr}(\theta \in A\mid \theta \in D)}{\mathsf{Pr}(\theta \in B \mid \theta \in D)} \leq n^{-1} \} \in \mathcal{U}.$$ The aim is to build an ultrafilter $\mathcal{U}$ for which this holds. 

For any  $n \in \mathbb{N}$, define $$C^n_{AB} \equiv \{ D \in [V]^{< \omega} :  \frac{\mathsf{Pr}(\theta \in A\mid \theta \in D)}{\mathsf{Pr}(\theta \in B \mid \theta \in D)}  \leq n^{-1}  \}.$$

Moreover, let  $$ A_x \equiv \{ D \in [V]^{< \omega} :  x \in D \}  .$$

Define also $$\mathcal{F} \equiv \{ C^n_{AB} : n\in \mathbb{N},  \left \vert A \right \vert < \left \vert B \right \vert \} \cup \{A_x : x \in V \}.$$
 We want to prove that  $\mathcal{F}$ has the finite intersection property. Therefore take any $x_1,\ldots ,x_k \in V$, and any $\langle A_1,B_1, n_1 \rangle, \ldots, \langle A_l,B_l,n_l \rangle$ such that $\left \vert A_j \right \vert < \left \vert B_j  \right \vert $ and $n_j \in \mathbb{N}$ for $j \leq l.$ Assume for the construction that $|A_1|\leq|A_2|\leq\dots\leq|A_l|$. For every finite $D$, if $ \{ x_1,\ldots . x_k \} \subseteq D,$ then $D \in \bigcap_{i\leq k} A_{x_i} .$ So setting $n=max\{n_j:j<l\}$ we will extend $\{ x_1,\ldots . x_k \}$ to a set in $C^n_{A_jB_j}$ , and hence $C^{n_j}_{A_jB_j}$ , for each $j \leq l$.
Set $F_0= \{ x_1,\ldots . x_k \}$ and $a_0=|F_0\cap A_1|$. As $B_1$ is infinite and of larger cardinality than $A_1$ we add $n \cdot a_0$ elements of $B_1 \setminus A_1$ to $ F_0$, yielding a finite set $F_1$. Now set $a_1=|F_1\cap A_2|$, and add $n \cdot a_1$ elements of $B _2\setminus (A_1\cup A_2)$  to $F_1$ to give $F_2$. Note we can find these elements of $B_2$ as $|B_2|>|A_2|\geq|A_1|$. Continuing in this manner, set $F=F_l$.
Then we have ensured that for all $j\leq l$ $$\frac{\mathsf{Pr}(\theta \in A_j\mid \theta \in F)}{\mathsf{Pr}(\theta \in B_j \mid \theta \in F)}  \leq n^{-1}, $$ and so we have $F \in C^n_{A_jB_j},$ and since $D \subseteq F,$ we also have $F \in \bigcap_{i\leq k} A_{x_i} .$ 

So $\mathcal{F}$ indeed has the finite intersection property, whereby it can be extended to a filter and then further to an ultrafilter $\mathcal{U}$. By design, then, the resulting probability function $\mathsf{Pr}_{\mathcal{U}}$ is super-regular.
\end{proof}
\end{theorem}

Once again, Hume's Principle for probability cannot hold for the notion of probability that we are investigating. But this leaves open the question whether the \emph{converse} of Hume's Principle for probability can be made to hold. This is called  \emph{Cantor's Principle} in \cite{BDF}, where the authors investigate it in the context of their Euclidean theory of size:
\begin{definition}[Cantor's Principle]
 $$ \mathsf{Pr}_{\mathcal{U}}(\theta \in A) = \mathsf{Pr}_{\mathcal{U}}(\theta \in B)  \Rightarrow \left \vert A \right \vert = \left \vert B \right \vert   .$$
 \end{definition}
\noindent Benci, Forti, and Di Nasso prove that `Cantor's Principle' can be made to hold \cite[section 3.2]{BDF}. It is also clear that Cantor's Principle follows from super-regularity.


\subsubsection{The power set principle}

The question whether
$$  \forall A, B\in V: \vert A \rvert < \lvert B \rvert \Rightarrow   \vert \mathcal{P}(A) \rvert < \lvert \mathcal{P}(B) \rvert  $$
\noindent is true, is independent of the axioms of set theory. (Of course the principle is true if the Generalised Continuum Hypothesis holds.)
Like the cardinality operator, our $\textrm{NAP}$ probability functions are measures of some kind. One might wonder what should  follow from $  \mathsf{Pr}_{\mathcal{U} }(\theta \in A) < \mathsf{Pr}_{\mathcal{U}}( \theta \in B ). $ In particular, given that $\mathsf{Pr}_{\mathcal{U} }$ is intended to be a \emph{fine-grained} quantitative possibility measure, perhaps probability should be expected to co-vary with the power set operation in some fairly direct manner. In other words, it is natural to ask if the following principle can be made to hold:
\begin{definition}[Power Set Condition]\label{power set condition}\index{power set condition|textbf}
$$ \forall A,B \in V:  \mathsf{Pr}_{\mathcal{U} } (\theta \in A) < \mathsf{Pr}_{\mathcal{U} } (\theta \in  B ) \Leftrightarrow  \mathsf{Pr}_{\mathcal{U} } (\theta \in\mathcal{P}(A))   < \mathsf{Pr}_{\mathcal{U} }(\theta \in \mathcal{P}(B)). $$
\end{definition}

It turns out that the power set condition can indeed be satisfied:
\begin{theorem}\label{powersetcondition}
There are fine ultrafilters\index{ultrafilter} $\mathcal{U}$ such that $\mathsf{Pr}_{\mathcal{U}}$ satisfies the power set condition\index{power set condition}.
\end{theorem}

\noindent The argument for this is somewhat more involved.

We aim to prove Theorem \ref{powersetcondition} by building the probability function up from an ultrafilter $\mathcal{U}$ which is based on a pre-filter $\mathcal{C}\subseteq \mathcal{P}([V]^{<\omega})$ that has the finite intersection property. 

The class $\mathcal{C}$ is built up in stages, and in such a way that it eventually witnesses the truth of the power set condition for all $A,B \in V$.

\vspace{0.3cm}

\noindent \texttt{Stage 0}

\noindent The class $\mathcal{C}_0$ consists of all $$  A_x \equiv \{ a\in [V]^{<\omega} : x \in a\},$$ for $x \in V$. This is to ensure that the ultrafilter that will be built from $\mathcal{C}$ is fine. We know that $\mathcal{C}_0$ has the finite intersection property.

\vspace{0.3cm}

\noindent \texttt{Limit stages}

\noindent For limit stages $\lambda$, we simply set $\mathcal{C}_{\lambda} \equiv \bigcup_{\beta < \lambda}\mathcal{C}_{\beta}$.

\vspace{0.3cm}

\noindent \texttt{Successor stages}

\noindent Given fine-ness, we may, and will, ignore the elements of $V_{\omega}$.
 At stage $\alpha > \omega$, where $\alpha$ is a successor ordinal, we consider the sets of $V_{\alpha } \backslash V_{\alpha - 1}$ and ensure that the power set condition eventually holds for all these sets and their power sets, by adding families of finite sets to $\mathcal{C}_{\alpha - 1}$ in such a way that the finite intersection property is preserved. 

As an illustrative and indeed representative example we do the case where $\alpha = \omega + 1$.

Let there be given an enumeration $\{ A_1, B_1 \},\hdots, \{ A_{\beta}, B_{\beta} \}, \hdots$ of the pairs of elements of $V_{\omega + 1} \backslash V_{\omega}$.

For the induction, we assume that, by having added appropriate sets of finite sets to $\mathcal{C}_0$, the power set condition holds for $\{ A_1, B_1 \},\hdots, \{ A_{\beta}, B_{\beta} \}$ and their power sets, and that in the process the finite intersection property has been preserved. The aim is now to extend this so that it also holds for $\{ A_{\beta + 1}, B_{\beta + 1} \}$. 
In other words, we have constructed $\mathcal{C}_1^{\beta}$, and we want to obtain 
$\mathcal{C}_1^{\beta + 1}$, where $\mathcal{C}_1^0 \equiv \mathcal{C}_0$.

\begin{definition}
$$ C_{A<B} \equiv \{ D \in [V]^{< \omega}: \frac{\left \vert A \cap D  \right \vert } {\left \vert D  \right \vert } < \frac{\left \vert B \cap D  \right \vert } {\left \vert D  \right \vert } \}
  .$$ 
\end{definition}

\begin{definition}
$$ C_{A\geq B} \equiv \{ D \in [V]^{< \omega}: \frac{\left \vert A \cap D  \right \vert } {\left \vert D  \right \vert } \geq \frac{\left \vert B \cap D  \right \vert } {\left \vert D  \right \vert } \}
  .$$ 
\end{definition}

\vspace{0.2cm}

\noindent \texttt{Claim}

\noindent Either $\mathcal{C}_1^{\beta} \cup \{ C_{A_{\beta} < B_{\beta}} \}$ has the finite intersection property, or $\mathcal{C}_1^{\beta} \cup \{ C_{A_{\beta} \geq B_{\beta}} \}$ has the finite intersection property (or both).

\vspace{0.2cm}

\noindent \texttt{Proof}

\noindent Suppose not. Then there is a finite intersection $F$ of elements of $\mathcal{C}_1^{\beta}$ such that $F \cap C_{A_{\beta} < A_{\beta}}  = \emptyset$, and there is a finite intersection $F'$ of elements of $\mathcal{C}_1^{\beta}$ such that $F' \cap C_{A_{\beta} \geq B_{\beta}}  = \emptyset$. But then $( F \cap F') \cap C_{A_{\beta} < B_{\beta}}  = \emptyset$ and $(F \cap F') \cap C_{A_{\beta} \geq B_{\beta}}  = \emptyset$. But $C_{A_{\beta} < B_{\beta}} \cup  C_{A_{\beta} \geq B_{\beta}} = [V]^{< \omega}.$ So then $(F \cap F') = \emptyset$. But this contradicts the inductive assumption that $\mathcal{C}_1^{\beta}$ has the finite intersection property.

\vspace{0.3cm}

Thus define $\mathcal{C}^{\beta+1}_1$ to be $\mathcal{C}_1^{\beta} \cup \{ C_{A_{\beta} < B_{\beta}} \}$ if this has the finite intersection property, or $\mathcal{C}_1^{\beta} \cup \{ C_{A_{\beta} \geq B_{\beta}} \}$ otherwise, and by the claim, $\mathcal{C}^{\beta+1}_1$ has the finite intersection property.
Now setting $\mathcal{C}^-_1 \equiv \bigcup_{\beta}\mathcal{C}^{\beta}_1, $ we may conclude that $\mathcal{C}^-_1$ has the finite intersection property.

At this point we must extend $\mathcal{C}^-_1$ by adding to  $\mathcal{C}^-_1$:
\begin{itemize}
\item every set of the form $C_{\mathcal{P}(A)<\mathcal{P}(B)}$ such that $C_{A<B} \in \mathcal{C}^-_1$;
\item every set of the form $C_{\mathcal{P}(A) \geq \mathcal{P}(B)}$ such that $C_{A\geq B} \in \mathcal{C}^-_1$.
\end{itemize}
\noindent Call the resulting set $\mathcal{C}_1$. Our aim is to prove that $\mathcal{C}_1$ has the finite intersection property.

Consider an arbitrary non-empty finite family $\mathcal{F} \subseteq \mathcal{C}_1$. Without loss of generality we may assume that the `judgements' in  $\mathcal{F}$ of the form $C_{\mathcal{P}(A)<\mathcal{P}(B)}$ or $C_{\mathcal{P}(A) \geq \mathcal{P}(B)}$, taken together, describe a finite total pre-ordering relation $R$ on some set $\{ \mathcal{P}(A_1),\hdots, \mathcal{P}(A_k) \}$. Further, we may also assume that for and sets $A$ and $B$ from $V_{\omega + 1} \backslash V_{\omega}$,   $C_{\mathcal{P}(A)<\mathcal{P}(B)}\in\mathcal{F}$ if and only if $C_{A<B}\in\mathcal{F}$, and $C_{\mathcal{P}(A) \geq \mathcal{P}(B)}$ iff  $C_{A\geq B}\in\mathcal{F}$. Thus $\mathcal{F}$ contains witnesses for all the relevant judgements we may be interested in.

Let $\mathcal{F}^- = \mathcal{F} \cap \mathcal{C}^-_1$, so $\mathcal{F}^- $ consists only of judgements about sets in $V_{\omega + 1} \backslash V_{\omega}$. Then we know from the foregoing that $\bigcap \mathcal{F}^- \neq \emptyset$. 
So take some $F^- \in \bigcap \mathcal{F}^-$.
Our plan is inductively to extend $F^-$, using the pre-order $R$, to a finite set $F \in \bigcap \mathcal{F}$.

We will add to $F^-$ elements that ensure that the constraints of $R$ are satisfied. Moreover, by choosing the elements to be added to $F^-$ from $V_{\omega + 1} \backslash V_{\omega}$,\footnote{For later stages we will take these sets from $V_{\alpha+ 1} \backslash V_{\alpha}$, i.e. sets of rank $\alpha$.}  we ensure that the constraints imposed by $\mathcal{F}^-$ remain satisfied. As a result, $F$ will satisfy all constraints from $\mathcal{F}$, so $\bigcap \mathcal{F}\neq \emptyset$ and hence $\mathcal{C}_1$ has the finite intersection property.

As an example, suppose that 
$$ R = \mathcal{P}(A_1) <  \mathcal{P}(A_2) < \mathcal{P}(A_3) = \mathcal{P} (A_4)   .$$

\vspace{0.2cm}

\noindent (1) We start by ensuring that $\mathcal{P}(A_1) <  \mathcal{P}(A_2)$ is satisfied. 

Suppose that $F^-$ already contains $n$ elements of $\mathcal{P}(A_1)$. Since $C_{A_1 < A_2} \in \mathcal{F}$, there must be an element $x^-\in A_2 \backslash A_1$. This implies that there are infinitely many infinite sets $x$ in $\mathcal{P}(A_2) \backslash \mathcal{P}(A_1)$ such that $x^- \in x$: we add $n+1$ such elements to $F^-$, and call the resulting finite set $F^-_1$.

\vspace{0.2cm}

\noindent (2) We proceed in  similar fashion to ensure that $\mathcal{P}(A_2) <  \mathcal{P}(A_3)$ is satisfied:

Suppose that $F^-_1$ already contains $m$ elements from  $\mathcal{P}(A_2)$, observing that it may be the case that $m> n+1$, for there may already be a finite number of elements of $\mathcal{P}(A_2)$ in $F^-$. Since $C_{A_2 < A_3} \in \mathcal{F}$, there must be an element $y^-_1\in A_3 \backslash A_2$, and since $C_{A_1 < A_3} \in \mathcal{F}$, there must be an element $y^-_2\in A_3 \backslash A_1$. So there are infinitely many infinite sets $y$ in $\mathcal{P}(A_3)$ such that $y^-_1, y^-_2 \in y$: add $m+1$ such elements to $F^-_1$, and call the resulting set $F^-_2$.

\vspace{0.2cm}

\noindent (3)
Now suppose that there are $m_1$ elements of $\mathcal{P}(A_3)$ in $F^-_2$, and $m_2$ elements of $\mathcal{P}(A_4)$ in $F^-_2$. Moreover, suppose that $m_2 < m_1$. (The case where $m_1 < m_2$ is similar.) Since $C_{A_3 \geq A_4},C_{A_4 \geq A_3}\in \mathcal{F}$, but also $A_3 \neq A_4$, there must be some $x_1 \in A_3 \backslash A_4$ and some $x_2 \in A_4 \backslash A_3$. Moreover, since $C_{A_1 < A_4}, C_{A_2 < A_4} \in \mathcal{F}$, there are elements $x_3 \in A_4 \backslash A_1, x_4 \in A_4 \backslash A_2$. So $\mathcal{P}(A_4)$ contains infinitely many infinite sets $x$ such that $\{ x_2,x_3,x_4  \}\subset x$. Similarly, $\mathcal{P}(A_3)$ contains infinitely many infinite sets $x$ that are outside $\mathcal{P}(A_1), \mathcal{P}(A_2), \mathcal{P}(A_4)$. So we add a sufficient number of such elements to $F^-_2$ so that there are an equal number $p$ of ``witnesses'' for $\mathcal{P}(A_3)$ as for $\mathcal{P}(A_4)$ but where $p$ is larger than the number of witnesses for $\mathcal{P}(A_2)$. Call the resulting set $F^-_3$.

\vspace{0.2cm}

\noindent (4) To conclude, we set $F \equiv F_3^-$. It is clear that $F \in \bigcap \mathcal{F}$. 

\vspace{0.2cm}

This procedure of extending $F^-$ easily generalises to any finite total pre-ordering on $\{ \mathcal{P}(A_1),\hdots, \mathcal{P}(A_k) \}$. Thus we have shown that $\mathcal{C}_1$ has the finite intersection property.

\vspace{0.3cm}

This procedure for extending $\mathcal{C}_0$ to $\mathcal{C}_1$ while preserving the finite intersection property also works for larger successor ordinals: at level $V_{\alpha+1}$ (stage $\beta+1$ with $\alpha=\omega+\beta$)  we can extend the corresponding $F^-$ using subsets of rank $\alpha$. As we have said above, at limit stages we can simply take unions. Ultimately we set $\mathcal{C} \equiv \bigcup_{\alpha \in On}\mathcal{C}_{\alpha}$.

The class $\mathcal{C}$ will then have the finite intersection property, so it can be extended to a filter and then to an ultrafilter $\mathcal{U}$. The probability function based on $\mathcal{U}$ will make the power set condition true for all $A,B \in V$, and this concludes the proof of theorem \ref{powersetcondition}.

\vspace{0.4cm}

Our proof actually shows something slightly stronger: for all $A,B$ with $\left \vert A \right \vert , \left \vert B \right \vert \geq \omega$, we have $$ \mathsf{Pr}_{\mathcal{U} } (\theta \in A) < \mathsf{Pr}_{\mathcal{U} } (\theta \in  B ) \Leftrightarrow  \mathsf{Pr}_{\mathcal{U} } (\theta \in\mathcal{P}(A))  \ll \mathsf{Pr}_{\mathcal{U} }(\theta \in \mathcal{P}(B)) .$$ The reason is that in enlarging the set $F^-$ we always have infinitely many elements to choose from.

For any probability measure $\mathsf{Pr}_{\mathcal{U}}$ that satisfies power set condition we also have that $ \forall A,B \in V, \forall n\in\omega$:
$$\mathsf{Pr}_{\mathcal{U}} (\theta \in A) < \mathsf{Pr}_{\mathcal{U}}(\theta \in B ) \Leftrightarrow  \mathsf{Pr}_{\mathcal{U}}(\theta \in \mathcal{P}^{n}(A))   < \mathsf{Pr}_{\mathcal{U}}(\theta \in\mathcal{P}^{n}(B)) $$
where $\mathcal{P}^{n}(A)=\mathcal{P}(\mathcal{P}(\dots\mathcal{P}(A)\dots))$. An easy argument shows this cannot extend to infinite applications of the power set operation.

One might wonder whether the motivations behind the power set condition should not also support imposing the following \emph{restricted power set condition} on $\mathsf{Pr}_{\mathcal{U}}$:\footnote{Thanks to Philip Welch for this question.}
\begin{question}
Are there probability measures such that $$ \forall A,B \in V:  \mathsf{Pr}_{\mathcal{U}}(\theta \in A) < \mathsf{Pr}_{\mathcal{U}}(\theta \in B ) \Leftrightarrow  \mathsf{Pr}_{\mathcal{U}}(\theta \in  [A]^{< \omega})   <\mathsf{Pr}_{\mathcal{U}}(\theta \in [B]^{< \omega})?$$
\end{question}


\subsection{The ordinals}

For $\alpha \geq \omega , $ in each level $V_{\alpha + 1} \setminus V_{\alpha}$ of the iterative hierarchy one finds only one ordinal, but infinitely many sets that are not ordinals. This might lead one to believe that a probability function on $V$ should satisfy $$ \mathsf{Pr}_{\mathcal{U}}(\theta \in \mathrm{On}) \approx 0 , $$ where `On' is the class of ordinals.

Just as it seems reasonable to require that the probability of choosing an even natural number from the set of natural numbers must be equal to or infinitesimally close to $\frac{1}{2}$ (see \cite[section 6.2]{WH}), it seems reasonable to require that
$$  \mathsf{Pr}_{\mathcal{U}}(\theta \in \mathrm{Even} \mid \theta \in \mathrm{On}) \approx \frac{1}{2} , $$ where `Even' is the class of even ordinals, which is defined in the obvious way.

Moreover, between any two limit ordinals there are infinitely many successor ordinals, so one might expect
$$  \mathsf{Pr}_{\mathcal{U}}(\theta \in \mathrm{Lim} \mid \theta \in \mathrm{On}) \approx 0 , $$
where `Lim' is the class of limit ordinals.

We will sketch how probability functions can be constructed that meet these expectations. Indeed, we will see that there are probability functions that meet these `ordinal expectations' and in addition meet the size constraint of super-regularity.

\begin{theorem}
There are super-regular probability functions $Pr$ such that:
\begin{enumerate}
\item $ \mathsf{Pr}_{\mathcal{U}}(\theta \in \mathrm{On}) \approx 0;$
\item $ \mathsf{Pr}_{\mathcal{U}}(\theta \in \mathrm{Even} \mid \theta \in \mathrm{On}) \approx 2^{-1} ;$
\item $ \mathsf{Pr}_{\mathcal{U}}(\theta \in \mathrm{Lim} \mid \theta \in \mathrm{On}) \approx 0. $
\end{enumerate}
\begin{proof}
As before, the aim is wisely to choose the ultrafilter $\mathcal{U}$ on which $\mathsf{Pr}_{\mathcal{U}}$ is based. We want $\mathcal{U}$ to be such that for all $k,l,m \in \mathbb{N}$:
\begin{itemize}
\item $\frac{\mathsf{Pr}_{\mathcal{U}}(\theta \in A)}{\mathsf{Pr}_{\mathcal{U}}(\theta \in B)} \leq k^{-1}$ if $\omega \leq \left \vert A \right \vert < \left \vert B \right \vert ;$
\item $   \mathsf{Pr}_{\mathcal{U}}(\theta \in \mathrm{Even} \mid \theta \in \mathrm{On}) - \mathsf{Pr}_{\mathcal{U}}(\theta \in \mathrm{Odd} \mid \theta \in \mathrm{On}) \leq l^{-1} $ and $\mathsf{Pr}_{\mathcal{U}}(\theta \in \mathrm{Lim} \mid \theta \in \mathrm{On}) \leq  l^{-1}  ;$
\item $  \mathsf{Pr}_{\mathcal{U}}(\theta \in On) \leq m^{-1} .$
\end{itemize}
Now we define:
\begin{itemize}
\item $ A_x \equiv \{ D \in [V]^{< \omega} :  x \in D \} ;$
\item $  C^k_{AB} \equiv \{ D \in [V]^{< \omega} :  \frac{\mathsf{Pr}[A\mid D]}{\mathsf{Pr}[B \mid D]}  \leq k^{-1}  \}  ;$
\item $ I^l \equiv \{ D \in [V]^{< \omega} :  \forall \alpha \in D \exists \beta \exists n \geq l ( \alpha \in [ \beta, \beta + n ] \subseteq D )  \} ;$
\item $W^m \equiv \{ D \in [V]^{< \omega} : \mathsf{Pr}[On \mid D] \leq m^{-1} \} .$
\end{itemize}
And now we set: $$  \mathcal{F}_0 \equiv \{A_x, C^k_{AB}, I^l, W^m: x \in V, k,l,m \in \mathbb{N} \textrm{ and } \omega \leq \left \vert A \right \vert < \left \vert B \right \vert \}    $$
\textbf{Claim:} $  \mathcal{F}_0$ has the finite intersection property. 

Let some $x_1,\ldots , x_n $ be given.  Now  $\bigcap_{i\leq n} I^{l_i}=I^l$ where $l={max\{l_i:i<n\}}$, and similarly for $\bigcap_{i\leq n} W^{m_i}$, so as before in theorem \ref{superregtheorem}, it suffices to concentrate on the highest values of $k,l,m$.

\noindent (1) $A \in \bigcap_{i\leq n}A_{x_i} \Leftrightarrow \{ x_1,\ldots , x_n\} \subseteq A .$ 
So we start with the finite set $A_0 \equiv  \{ x_1,\ldots , x_n \} ,$ and will extend it. 

\noindent (2) Again we concentrate on one pair $\langle A, B\rangle$ such that $\omega \leq \left \vert A \right \vert < \left \vert B \right \vert$; we leave out further cases as they are similar. There are arbitrarily large finite subsets $C \subseteq B$ that are $l$-isolated from elements of $A_0$, meaning that each ordinal in $C$ is more than $l$ ordinals removed from any ordinal in $A$. We choose any such $C \subseteq B$ that is of size at least $k \cdot n$, and we set $A_1 \equiv A_0 \cup C$.

\noindent (3) Now we extend $A_1$ to ensure that all ordinal intervals are of length $\geq l$: for each $\alpha \in A_1$, we add $ \alpha +1,\ldots , \alpha +l  $. Call the resulting finite collection $A_2$.
Note that by our choice of $l$-isolated elements in (2), none of $ \alpha +1,\ldots , \alpha +l  $ are elements of $A$.

\noindent (4) Let $\left \vert A_2 \right \vert = j$. Then we add $j \cdot m $ elements of $V \setminus ( A \cup B \cup On   )$ to $A_2$ and call the resulting set $A_3$.

It is now routine to verify that $A_3 \in \bigcap_{i\leq n}A_{x_i} \cap C^k_{AB} \cap I^l \cap W^m$. The case including further sets $C^k_{A'B'}$ is  similar, thus the claim is verified. So $\mathcal{F}_0$ indeed has the finite intersection property, whereby it can be extended to a filter and then further to an ultrafilter $\mathcal{U}$. By design, the resulting probability function $\mathsf{Pr}_{\mathcal{U}}$ has the required properties.
\end{proof}
\end{theorem}


\section{The bootstrapping approach}\label{bootstrapping}

The probability $\mathsf{Pr}_{\mathcal{U}}(\theta \in A)$ is obtained by `summing up' the probabilities $\mathsf{Pr}(\theta \in A \mid \theta \in S)$ for all `small' parts $S$ of $V$; such  $\mathsf{Pr}(\theta \in A \mid \theta \in S)$ are seen as \emph{approximations} of $\mathsf{Pr}_{\mathcal{U}}(\theta \in A)$.

In the finite snapshot approach, `small' in this context means `finite'.
But from a conceptual point of view, `finite' might be taken to be \emph{too} small as far as the test sets (or snapshots) are concerned.
Compared to $V$, \emph{all} sets ---and not just the finite sets--- are small. So to determine $\mathsf{Pr}_{\mathcal{U}}(\theta \in A)$, we should take the `limit' of the values $\mathsf{Pr}(\theta \in A \mid \theta \in S)$, where $S$ is a set of any size. Then if $S$ is infinite, $\mathsf{Pr}(\theta \in A \mid \theta \in S)$ cannot just be taken to be given by the ratio formula but needs to be \emph{defined}.

In the approach to which we now turn (the bootstrapping approach), a probability $\mathsf{Pr}_{\mathcal{U}}(\theta \in A)$ is determined by the probabilities $\mathsf{Pr}_{\mathcal{U}}(\theta \in A \mid \theta \in S)$, where $\mathsf{Pr}_{\mathcal{U}}(\theta \in A \mid \theta \in S)$, for $S$ a large set, is then in turn determined by probabilities $\mathsf{Pr}_{\mathcal{U}}(\theta \in A \mid \theta \in S')$ for $S'$ being smaller `snapshots' than $S$, and so on, until we reach the finite snapshots and can appeal to the probability functions that were discussed in the previous sections. Thus the bootstrapping account can be seen as a \emph{generalisation} of the finite snapshot approach.


\subsection{The rough idea}

In general terms, this is how we will proceed:

\noindent (1) By the construction from the previous section, a fine ultrafilter on $[S]^{<\omega}$ yields a notion of probability on all sets $S \in V$ with $\left \vert S \right \vert  < \omega_{1}$. In other words, this yields a suitable notion of probability, call it $\mathsf{Pr}^S$, for every \emph{countable} set $S$.

\noindent (2) The notion of $\mathsf{Pr}^S$ for all $S \in V$ with $\left \vert S \right \vert < \omega_2$ is determined using the notion of probability on \emph{countable} sets: the probability of $A$ on such an $S$ is determined by the class of probabilities of $A$ on the countable `snapshots' of $S$. Using these countable probability functions, a fine ultrafilter on $[S]^{< \omega_1}$ gives us a notion of probability on sets $S$ with $\left \vert S \right \vert < \omega_2$.

Again the resulting functions $\mathsf{Pr}^S$ are essentially NAP-functions as defined in \cite{BHW}. They are total, regular, etc.

\vspace{0.3cm}

\ldots

\vspace{0.3cm}

\noindent ($\beta$) A fine ultrafilter on $[S]^{< \omega_{\alpha}}$, together with probability functions $\mathsf{Pr}^S$ for all $S$ such that $\left \vert S \right \vert < \omega_\alpha$, yields a notion of probability on all sets $S$ with $\left \vert S \right \vert < \omega_{\alpha +1}$.

\vspace{0.3cm}

\ldots

\vspace{0.3cm}

Limit stages of course do not present a problem.
So by transfinite recursion on cardinality this yields for \emph{every} set $S$ a notion $\mathsf{Pr}^S$ of probability on $S$.

Then a fine ultrafilter $\mathcal{U}$ on $V=[V]^{<{Card}}$ yields, using the general notion $\mathsf{Pr}^S$ for $S\in V$, a notion $\mathsf{Pr}^V$ that is a total (class) function from properties $A$ and random variables $\theta$ to values $\mathsf{Pr}^V(\theta \in A)$ in  a non-Archimedean class field. This probability function again satisfies the principles of the theory NAP in \cite{BHW}.

For this construction, what we need is suitable (fine) ultrafilters on small, and somewhat larger, and large, $\hdots$ sets, and a fine ultrafilter $\mathcal{U}$ on $[V]^{<Card}$. But we will see that all the set ultrafilters used in the construction can be \emph{uniformly} obtained as \emph{restrictions} to sets $S$ of the given fine ultrafilter on $[V]^{<Card}$. So $\mathsf{Pr}^V$ is determined by \emph{one} initial choice of $\mathcal{U}$, whereby $\mathsf{Pr}^V$ can be seen as the `limit' of its set-restrictions $\mathsf{Pr}^S$, where the functions $\mathsf{Pr}^S$ can in turn be seen as `limits' of restrictions to \emph{their} small subsets. This uniform construction has the advantage that the resulting probability functions are all \emph{coherent}, in the sense that for a set $T$, $\mathsf{Pr}^S(A|T)$ is the same for all $S\supseteq T$ and hence also for $V$.

Now it is time to look at details of the construction.


\subsection{Details 1: Restrictions of fine ultrafilters}

Since our construction involves ultrafilters on sets $[S]^{< \kappa}$ with $\kappa > \omega$, we make the following definition, which accords with the usual definition of fineness on $[S]^{< \omega}$.

\begin{definition}
For any infinite cardinal $\kappa$, an ultrafilter on $[S]^{< \kappa}$ is fine iff for every $x \in S:$  $$ \{T \in [S]^{< \kappa}: x \in T \} \in \mathcal{U}.$$
The notion of `set-fine' ultrafilter on $V$ is defined in the obvious way.
\end{definition}

We first show that appropriate restrictions of ultrafilters to smaller sets can be obtained in a uniform fashion.

\begin{definition}
Suppose $S \in V$, $\lvert S \rvert = \kappa$, and $\mathcal{U}$ a fine ultrafilter on $[S]^{<\kappa}$, and $S' \subseteq S$ with $\lvert S' \rvert = \alpha < \kappa$.  Then we define the restriction $\mathcal{U}_{S'}$ of $\mathcal{U}$ to $S'$ as follows. 

For any $X \in \mathcal{P}([S]^{<\kappa})$, let
\[
X_{S'} \equiv \{ y \mid   \exists z \in X: y = z \cap S' \textrm{ and } \lvert y \rvert < \alpha \}.
\]
Then $\mathcal{U}_{S'} \equiv \{ X_{S'} \mid X \in \mathcal{U} \}.$ 

\end{definition}

\begin{proposition}\label{restrict}
For any $S \in V$ with $\left \vert S \right \vert = \kappa $, there are fine ultrafilters $\mathcal{U}$ on $[S]^{<\kappa}$ that restrict to a fine ultrafilter on every $S' \subseteq S$ with $\left \vert S' \right \vert = \alpha$, and $\omega \leq \alpha < \kappa$.

Further, such ultrafilters are \emph{coherent} in that if $T\subset S'$ with $\omega\leq |T| <|S'|$, then $(\mathcal{U}_{S'})_T=\mathcal{U}_T$.

\begin{proof}
We build the ultrafilter from a pre-filter $\mathcal{F}_0$  (i.e., a set with the finite intersection property), which can then be extended to a filter and then to an ultrafilter.

For each $x \in S$, let $$ A_x \equiv \{ X \in [S]^{< \kappa}: x \in X   \}  .$$ 

And let for each $S'$ with $\left \vert S' \right \vert = \alpha < \kappa$ and $S' \subseteq S$: 
$$R^{S'} \equiv \{ X \in [S]^{< \kappa}:  X\cap S'\in[S']^{< \alpha} \}   .$$ 
Now set $$ \mathcal{F}_0 \equiv  \{ A_x: x \in S  \} \cup \{ R^{S'}: S' \subseteq S \textrm{ and } \left \vert S' \right \vert < \kappa  \} .$$ It is easy to see that $ \mathcal{F}_0$ has the finite intersection property and so can be extended to an ultrafilter $\mathcal{U}$. And by design, $\mathcal{U}$ is fine.

Clearly $\mathcal{U}_{S'} \subseteq \mathcal{P}([S']^{<\alpha}).$ We must check the fine ultrafilter properties:

\noindent (1) \textbf{Fine.} This follows from the fact that $\mathcal{U}$ is fine: for $x\in S'$ this is witnessed by $(A_x)_{S'}$.

\noindent (2)  \textbf{Finite intersection.} Let $X,Y \in \mathcal{U}_{S'}$. Then there are  $\overline{X}, \overline{Y} \in \mathcal{U}$ such that $X = \overline{X}_{ S'}$ and  $Y = \overline{Y}_{ S'}$. By the finite intersection property of $\mathcal{U}$, we know that $\overline{X} \cap \overline{Y} \in \mathcal{U}.$ But $X \cap Y \supseteq(\overline{X} \cap \overline{Y} )_{ S'}.$ So $X \cap Y \in \mathcal{U}_{S'}$.

\noindent(3)  \textbf{Ultra.} Take any $X \subseteq [S']^{< \alpha}$, and let $X^c \equiv [S']^{< \alpha} \backslash X.  $ 
Let $\overline{X} \equiv \{ x \in [S]^{< \kappa} \mid x \cap S' \in X  \}$ and let $\overline{X^c} \equiv \{ x \in [S]^{< \kappa} \mid x \cap S' \not \in X \}. $ 
Then $\overline{X^c} = [S]^{< \kappa} \backslash \overline{X}.$ By the ultra property for $\mathcal{U}$, we have $\overline{X} \in \mathcal{U}$ or $\overline{X^c} \in \mathcal{U}$. But $X = \overline{X}_{S'}$ and $X^c = \overline{X^c}_{S'}$. So $X \in \mathcal{U}_{S'}$ or $X^c \in \mathcal{U}_{S'}.$

\noindent (4)  \textbf{Non-principality.} This is implied by fineness.

\noindent (5)  \textbf{Empty set property}: We have to show that $\emptyset \not \in \mathcal{U}_{S'}$. It suffices to show that for each $X \in \mathcal{U}$, $X_{S'} \neq \emptyset$. Since $R^{S'} \in \mathcal{U}$, $X \cap R^{S'} \neq \emptyset$. But for any set $x$ in this intersection, $x\cap S'\in[S']^{< \alpha}$. So $x\cap S'\in X_{S'}\neq \emptyset.$

\vspace{0.2cm}

For coherence, take $T\subset S'\subset S$ with $|T|<|S'|<|S|$ and let $X\in \mathcal{U}$. As $R^{S'}\in\mathcal{U}$ it is enough to show that $((X\cap R^{S'})_{S'})_T=(X\cap R^{S'})_T$. Now $((X\cap R^{S'})_{S'})_T=\{y \mid   \exists z \in X\cap R^{S'}: y = z \cap T, |y|<|T| \textrm{ and } \lvert z\cap S' \rvert < |S'|\}$, but by definition, for any $z\in R^{S'}$ we have $\lvert z\cap S' \rvert < |S'|$. Thus   $((X\cap R^{S'})_{S'})_T=\{y \mid   \exists z \in X\cap R^{S'}: y = z \cap T \textrm{ and } |y|<|T| \}=(X\cap R^{S'})_T$.
\end{proof}

\end{proposition}

But this means that this property must also hold for fine ultrafilters on $[V]^{< Card}:$

\begin{consequence}
There are fine ultrafilters $\mathcal{U}$ on $[V]^{< Card}$, such that for every set $S$  with $\left \vert S \right \vert = \alpha$,  $\mathcal{U}_{S}$ is a fine ultrafilter on $[S]^{<\alpha}$ and the coherence property holds.

\begin{proof}
By the same reasoning as in the previous proposition.
\end{proof}
\end{consequence}


\subsection{Details 2: defining probability functions}\label{def}

Now we show how for every set, a probability function on that set can be defined. The same procedure can then be used to define a probability function on $V$, and these probability functions are \emph{coherent}.

The key is to spell out what is involved in the $\beta$-th step of the recursive procedure for defining probabilities on sets:

\vspace{0.2cm}

\noindent ($\beta$) A fine ultrafilter $\mathcal{U}$ on $[S]^{< \omega_{\beta}}$ (with $\omega_{\beta} = \left \vert S \right \vert$), together with probability functions $\mathsf{Pr}^T$ for all $T$ such that $\left \vert T \right \vert < \omega_{\beta}$, yields a notion of probability $\mathsf{Pr}^S$ on $S$.
\vspace{0.2cm}

As in section \ref{snapshot}, we define a function $f_{\theta \in A}$ such that for all $T \in [S]^{< \omega_{\beta} }$: 
$$  f_{\theta \in A} (T) \equiv \mathsf{Pr}^T (\theta \in A\cap T)  .$$
Similarly, we define a function $f_{\theta \in A \wedge \nu \in B}$ such that for all $T \in [S]^{< \omega_{\beta} }$:
$$  f_{\theta \in A \wedge \nu \in B} (T) \equiv \mathsf{Pr}^T (\theta \in A\cap T \wedge \nu \in B\cap T)  .$$

\noindent Then $\mathsf{Pr}^S (\theta \in A)$ is defined as $[f_{\theta \in A}]_{\mathcal{U}}$, and $ \mathsf{Pr}^S (\theta \in A \mid \nu \in B) $ is defined as $$\frac{[f_{\theta \in A \wedge \nu \in B}]_{\mathcal{U}}}{[f_{\nu \in B}]_{\mathcal{U}}}.$$ This function $\mathsf{Pr}^S$ will then be an NAP probability function in the sense of \cite{BHW}.



Now in an exactly similar way, we define a class probability function $\mathsf{Pr}^+_{\mathcal{U}}$ on $V$, using the probability functions on `small' classes (i.e., sets) and ultrafilters on `small' classes which (given proposition \ref{restrict}) we can now assume to have been defined on the basis of an ultrafilter $\mathcal{U}$ on $[V]^{<Card}$ with which we start.
The function $\mathsf{Pr}^+_{\mathcal{U}}$ is total, regular, and uniform for the same reasons as why its `smaller cousin' $\mathsf{Pr}_{\mathcal{U}}$ has these properties. 

We now check coherence. We will do this only for straight probabilities rather than random variables in general, as although coherence holds for random variables also, it is much more technical to state. Below we use $\mathsf{Pr}(A)$ to denote $\mathsf{Pr}(\iota\in A)$ where $\iota$ s the identity random variable.

\begin{proposition}
For any class $A$ and sets $T\subset S$ with $|T|<|S|$ we have $$\mathsf{Pr}^T(A)=\mathsf{Pr}^S(A|T).$$

\begin{proof}
We show by induction on $|T|$ that that the above holds for all $S\supset T$ with $|S|>|T|$. Strictly speaking, the range of $\mathsf{Pr}^T$ may be a different non-archimedean field to the range of $\mathsf{Pr}^S$, but there is a natural embedding of the former into the latter defined by $i([f]_{\mathcal{U}_T})=[\bar{f}]_{\mathcal{U}_S}$ where for $X\in S^{<|S|}$, $\bar{f}(X)=f(X\cap T)$. This is well-defined as $\{X\in S^{<|S|}:|X\cap T|<|T|\}=(R^T)_S\in\mathcal{U}_S$.

Using this embedding we have $i(\mathsf{Pr}^T(A))=i([f_A]_{\mathcal{U}_T})=[\bar{f_{A}}]_{\mathcal{U}_S}$. Now for $X\in(R^T)_S(\in{\mathcal{U}_S})$ we have:
$$\bar{f_{A}}(X)=f_A(X\cap T)=f_{A\cap T}(X\cap T)=\mathsf{Pr}^{X\cap T}(A\cap T).$$ As $X\in(R^T)_S$ we have $|X\cap T|<|T|$ so by our inductive hypothesis $$\mathsf{Pr}^{X\cap T}(A\cap T)=\mathsf{Pr}^X(A\cap T|T)=\frac{f_{A \cap T}(X)}{f_{T}(X)}.$$ But by definition, $\big[\frac{f_{A \cap T}}{f_{T}}\big]_{\mathcal{U}_S}=\mathsf{Pr}^S(A|T)$, so $[\bar{f_{A}}]_{\mathcal{U}_S}=\mathsf{Pr}^S(A|T)$ and we're done.
\end{proof}

\end{proposition}



\subsection{Comparison of the finite snapshot approach and the bootstrapping approach}

In our definition of the probability of a set theoretic property, the probability $\mathsf{Pr}^+_{\mathcal{U}}(\theta \in A)$ of a property $A$ is determined by the probabilities $Pr^S(\theta \in A )$ of $A$ on large `snapshots' $S$, where a probability $Pr^S(\theta \in A )$ (for $S$ a large set) is then in turn determined by the probabilities $Pr^{S'}(\theta \in A$ for $S'$ being smaller `snapshots' than $S$, and so on. Conceptually, the definition in section \ref{def} is superior to the simpler definition suggested from section \ref{snapshot}: we want to take the behaviour of the property on as many and as large `snapshots` as possible into account. 

It is not straightforward to compare the simple and the more involved definition: the simple method is based on an ultrafilter on $[V]^{< \omega}$ whereas the more involved method is based on an ultrafilter on $V=[V]^{< Card}$.

The obvious suggestion is to base the comparison on the relation between a probability function determined by an ultrafilter $\mathcal{U}$ on $[V]^{< Card}$ and its \emph{restriction}\footnote{This is a different notion of restriction to that defined in the previous section as here we are only restricting the index, while the underlying class remains the same ($V$).}  to $[V]^{< \omega}$ defined as $\mathcal{U}\upharpoonright \omega=\{X\cap[V]^{< \omega}|X\in\mathcal{U}\}$. But:

\begin{proposition}
Not all ultrafilters on $[V]^{< Card}$ restrict to ultrafilters on to $[V]^{< \omega}$. 

\begin{proof}
Consider $\mathcal{A} \cup \overline{[V]^{< \omega}}$, where $\mathcal{A}$ is the set of atoms (guaranteeing fine-ness) and $\overline{[V]^{< \omega}}$ is the relative complement of ${[V]^{< \omega}}$ in ${[V]^{< Card}}$.  Then $\mathcal{A} \cup \overline{[V]^{< \omega}}$ has the finite intersection property and so can be extended to a fine ultrafilter $\mathcal{U}$ on ${[V]^{< Card}}$. But $\emptyset \in \mathcal{U}\upharpoonright \omega$. So $\mathcal{U}$ does not restrict to an ultrafilter on $[V]^{< \omega}$.
\end{proof}

\end{proposition}

On the other hand, every fine ultrafilter on $[V]^{< Card}$ restricting to an ultrafilter on $[V]^{< \omega}$ essentially \emph{is} an ultrafilter on $[V]^{< \omega}$:
\begin{proposition}
Suppose $\mathcal{U}$ is a fine ultrafilter on $[V]^{< Card}$ restricting to an ultrafilter $\mathcal{U} \upharpoonright \omega$ on $[V]^{< \omega}$. Then $[V]^{< \omega} \in \mathcal{U}$.

\begin{proof}
Since $\mathcal{U} $ is ultra, we have $[V]^{< \omega}\in \mathcal{U}$ or $\overline{[V]^{< \omega}}\in \mathcal{U}$. But if $\overline{[V]^{< \omega}}\in \mathcal{U}$, then $\emptyset \in \mathcal{U}\upharpoonright \omega$, so that $\mathcal{U}$ does not restrict, contradicting the assumption. So $[V]^{< \omega}\in \mathcal{U}$.
\end{proof}
\end{proposition}

This means that the \emph{essentially} involved probability functions on $V$ cannot be reduced to `simple' probability functions on $V$.


\section{Conclusion}

In this article we have explored two methods for modelling, by means of non-Archimedean probability functions, the properties of random variables ranging over the set theoretic universe: the finite snapshot method and the bootstrapping method. Concerning the finite snapshot method, we found that many of the probabilistic properties that seem intuitively plausible can be satisfied. The bootstrapping method is more satisfying from a conceptual point of view, but we have only been able to show that the resulting probability functions satisfy minimal requirements. So much work remains to be done.


\end{document}